\documentclass[reqno,centertags, 12pt, draft]{amsart}

\usepackage{amsmath,amsthm,amscd,amssymb}

\usepackage{latexsym}

\usepackage{mathrsfs}

\newcommand{\bbN}{{\mathbb{N}}}

\newcommand{\bbR}{{\mathbb{R}}}

\newcommand{\bbP}{{\mathbb{P}}}

\newcommand{\bbZ}{{\mathbb{Z}}}

\newcommand{\bbC}{{\mathbb{C}}}

\newcommand{\Jab}{{J(\{a(j)\},\{b(j)\})}}

\newcommand{\no}{\nonumber}

\newcommand{\ti}{\tilde  }

\newcommand{\ac}{\text{\rm{ac}}}

\newcommand{\beq}{\begin{equation}}

\newcommand{\eeq}{\end{equation}}

\newcommand{\ba}{\begin{align}}

\newcommand{\ea}{\end{align}}

\allowdisplaybreaks

\numberwithin{equation}{section}

\newtheorem{theorem}{Theorem}[section]

\newtheorem{proposition}[theorem]{Proposition}

\newtheorem{lemma}[theorem]{Lemma}

\newtheorem{corollary}[theorem]{Corollary}

\theoremstyle{definition}

\newtheorem{definition}[theorem]{Definition}

\theoremstyle{remark}

\newtheorem*{remark}{Remark}

\newtheorem*{remarks}{Remarks}

\begin{document}

\title[Sparse Trees with Finite Dimensions]{Singular Continuous and Dense 
Point Spectrum for Sparse Trees with Finite Dimensions}

\author[J.\ Breuer]{Jonathan Breuer*}

\thanks{$^*$ Institute of Mathematics, The Hebrew University of Jerusalem, Jerusalem, 91904, 
Israel. Email: jbreuer@math.huji.ac.il}

\date{May 15, 2006}

\dedicatory{Dedicated to S.~Molchanov on the occasion of his 65'th birthday}

\maketitle
\begin{abstract}
Sparse trees are trees with sparse branchings. The Laplacian on some of these trees can be shown to have 
singular spectral measures. We focus on a simple family of sparse trees for which the dimensions can be naturally defined 
and shown to be finite. Generically, this family has singular spectral measures and eigenvalues that are dense in 
some interval.
\end{abstract}
%%%%%%%%%%%%%%%%%%%%%%%%%%%%%%%%%%% INTRODUCTION %%%%%%%%%%%%%%%%%%%%%%%%%%%%%%%%%

\section{Introduction}
This paper extends and complements the paper \cite{breuer-CMP} in which the notion of sparse trees was introduced. Sparse 
trees are trees which have arbitrarily long `one-dimensional' segments (by which we mean - intervals of $\bbZ$), separated by 
occasional non-trivial branchings. It is shown in \cite{breuer-CMP} that, when these trees are spherically symmetric, one 
may decompose the Laplacian as a 
direct sum of Jacobi matrices which have sparse `bumps' off the diagonal. The spectral theory of these matrices 
is similar to that of one-dimensional Schr\"odinger operators with sparse potentials (see \cite{last-review} 
and references therein). In particular, matrices of this type 
exist for which the spectral measures are singular with respect to Lebesgue measure. 
These ideas make it possible to construct simple examples of trees for which the Laplacian has interesting spectral 
behavior. Several examples with singular continuous spectrum were presented in \cite{breuer-CMP}. 

In this paper we will be concerned with a family of sparse trees that `interpolates' between $\bbZ^+$ and the Bethe lattice. 
These trees can be obtained from the Bethe lattice by replacing an edge, at a distance $n$ from the root, by a segment of 
length $\sim \gamma^n$ for some fixed $\gamma>1$. 
While the Bethe lattice is infinite dimensional, a tree obtained in this manner can be shown to 
have dimensionality $=\frac{\log \gamma k}{\log \gamma}$, where $k$ is the connectivity of the original Bethe lattice. 
(For our definition of dimension see section 3). Thus, by 
letting $\gamma$ vary from $1$ to $\infty$, one gets a family of trees corresponding at one end 
($\gamma=1$) to the Bethe lattice, and at the other end ($\gamma=\infty$) to $\bbZ^+$. 

We shall analyze the spectral 
properties of the Laplacian on these trees with the help of the decomposition described above and some tools from the 
spectral theory of Schr\"odinger operators with sparse potentials. The constant branching, however, turns out to be a 
technical difficulty. We will bypass this difficulty by using an idea from \cite{zlatos} - namely, we shall 
impose a certain probability measure on these trees and prove an `almost sure' result.
 
It turns out that the situation for these finite dimensional structures is markedly different 
from the one for $\bbZ^d$. These trees (generically) have purely singular spectral measures and some dense point spectrum. 

In addition to the new result described above, we also use this opportunity to expand the discussion on the 
basic setting and on some of the examples presented in \cite{breuer-CMP}. Some basic facts that were briefly mentioned in 
that paper (such as the self-adjointness of the Laplacian on normal sparse trees), will be explained here in 
greater detail.

We remark that graphs with singular continuous \cite{graph-lap} and pure point spectrum \cite{teplyaev} are known to exist. 
In this context, the family of sparse trees is interesting in that, when varying two parameter sequences (namely - 
the branching size and the distances between branchings), one encounters a rich spectrum of phenomena. We note, in 
particular, the existence of examples with spectral measures of fractional Hausdorff dimensions (see \cite{breuer-CMP} and 
theorem \ref{finite-dimensional-trees} below). 

This paper is structured as follows. The next section presents the notion of sparse trees and the decomposition theorem 
that is basic for all that follows. Some simple results concerning spectral measures for the Laplacian on sparse trees are 
given in section 3. 
Section 4 describes the construction of the finite dimensional trees mentioned above and our results 
for them. As mentioned above, this paper uses some ideas and tools from the spectral theory of discrete one-dimensional 
Schr\"odinger operators. Relevant notions and results are presented in the appendix.

We are grateful to Michael Aizenman, Nir Avni, Vojkan Jak\v si\'c, Yoram Last, Barry Simon, Simone Warzel and Andrej Zlato\v s for 
useful discussions. We also wish to thank Michael Aizenman for the hospitality of Princeton, where some of this 
work was done. 

This research was supported in part by THE ISRAEL SCIENCE FOUNDATION (grant no. 188/02) and by Grant no. \mbox{2002068} 
from the United States-Israel Binational Science Foundation (BSF), Jerusalem, Israel.

%%%%%%%%%%%%%%%%%%%%%%%%%%%%%%%%%%%% SECTION 2 %%%%%%%%%%%%%%%%%%%%%%%%%%%%%%%%%%%%%%%%%%%%%%%%%%%%%%%%%%%%%%%%%%
\section{Sparse Trees}   

As noted in the introduction, basic to the analysis which follows is a certain decomposition of the Laplacian 
on a sparse tree. Since this is possible only when the tree has a certain spherical symmetry, we start with:
\begin{definition}[Spherically Homogeneous Rooted Tree]
A rooted tree is called spherically homogeneous (SH) (see \cite{bass}) if any vertex $v$ of generation $j$ is connected with 
$\kappa_j$ vertices of generation $j+1$. A locally finite spherically homogeneous tree is uniquely determined by the 
sequence $\{\kappa_j\}_{j=0}^\infty$. By locally finite we mean that the valence of any vertex is finite.
\end{definition} 

\setlength{\unitlength}{0.7cm}
\begin{figure}[!hbp]
\begin{picture}(5,5)(4,-3)
\put(3,0){\line(1,0){1}}
\put(3,0){\circle*{0.15}}
\put(4,0){\circle*{0.15}}
\put(4,0){\line(1,2){0.5}}
\put(4.5,1){\circle*{0.15}}
\put(4,0){\line(1,-2){0.5}}
\put(4.5,-1){\circle*{0.15}}
\put(4.5,1){\line(1,0){1}}
\put(5.5,1){\circle*{0.15}}
\put(5.5,1){\line(1,0){1}}
\put(6.5,1){\circle*{0.15}}
\put(6.5,1){\line(1,0){1}}
\put(7.5,1){\circle*{0.15}}
\put(4.5,-1){\line(1,0){1}}
\put(5.5,-1){\circle*{0.15}}
\put(5.5,-1){\line(1,0){1}}
\put(6.5,-1){\circle*{0.15}}
\put(6.5,-1){\line(1,0){1}}
\put(7.5,-1){\circle*{0.15}}
\put(7.5,-1){\line(1,1){0.6}}
\put(7.5,-1){\line(1,-1){0.6}}
\put(7.5,1){\line(1,1){0.6}}
\put(7.5,1){\line(1,-1){0.6}}
\put(8.1,1.6){\line(1,0){0.2}}
\put(8.5,1.6){\line(1,0){0.2}}
\put(8.9,1.6){\line(1,0){0.2}}
\put(8.1,0.4){\line(1,0){0.2}}
\put(8.5,0.4){\line(1,0){0.2}}
\put(8.9,0.4){\line(1,0){0.2}}
\put(8.1,-1.6){\line(1,0){0.2}}
\put(8.5,-1.6){\line(1,0){0.2}}
\put(8.9,-1.6){\line(1,0){0.2}}
\put(8.1,-0.4){\line(1,0){0.2}}
\put(8.5,-0.4){\line(1,0){0.2}}
\put(8.9,-0.4){\line(1,0){0.2}}
\put(8.1,1.6){\circle*{0.15}}
\put(8.1,0.4){\circle*{0.15}}
\put(8.1,-1.6){\circle*{0.15}}
\put(8.1,-0.4){\circle*{0.15}}
\end{picture}
\caption{An example of a SH rooted tree with $\kappa_0=1$, $\kappa_1=2$, $\kappa_2=\kappa_3=\kappa_4=1$, $\kappa_5=2$ 
\ldots}
\end{figure}
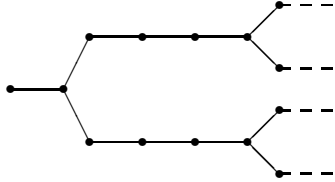

In the definition above, a vertex is said to be of generation $j$ if it is at a distance $j$ from the root - $O$
(where the distance between two vertices is defined as the number of edges of the unique path between them). Thus, for 
spherically homogeneous rooted trees, the valence of a vertex depends solely upon its location with respect to the root. 

Let $\{L_n\}_{n=1}^\infty$ and $\{k_n\}_{n=1}^\infty$ be two sequences of natural numbers such that $k_n \geq 2$ for all 
$n$, and $\{L_n\}_{n=1}^\infty$ is strictly increasing. We say that 
$\Gamma$ is a SH rooted tree of type $\{L_n,k_n\}_{n=1}^\infty$ if
\beq \label{kappa-j2}
\kappa_j= \left\{ \begin{array}{ll}
k_n & j=L_n \textrm{ for some } n\\
1 & \textrm{otherwise}
\end{array} \right.
\eeq
We say that $\Gamma$ is \emph{sparse} if $(L_{n+1}-L_n) \rightarrow \infty$ as $n \rightarrow \infty$. 
Since sparse trees are not regular (the coordination number is not constant), there are two natural choices for the 
Laplacian: 
\beq \label{lap1}
(\Delta f)(x)=\sum_{y : d(x,y)=1}f(y),
\eeq
and
\beq \label{lap2}
(\ti{\Delta} f)(x)=\sum_{y : d(x,y)=1 }f(y)-\#\{y : d(x,y)=1 \}\cdot f(x)
\eeq
where $\#A$ for a finite set $A$ is the number of elements in $A$ ($d(x,y)$ denotes the distance on the tree). 
For simplicity, we shall restrict our attention to $\Delta$, though we note that all our results hold 
(when properly modified) for $\ti{\Delta}$ as well.

It is clear that if $\{k_n\}_{n=1}^\infty$ is a bounded sequence then both $\Delta$ and $\ti{\Delta}$, 
on the tree, are bounded and self-adjoint. For unbounded coordination number, the issue of self-adjointness has to be 
addressed. 
\begin{definition} \label{normal-tree}
We call a SH rooted tree of type $\{L_n,k_n\}_{n=1}^\infty$ - $\Gamma$  - \emph{normal} if $\{k_n\}$ unbounded implies that 
$\limsup_{n \rightarrow \infty} (L_{n+1}-L_n)>1$. 
\end{definition}
The appendix has a proof that the Laplacians on normal SH  rooted trees are self-adjoint.
Clearly, any sparse tree is normal.

The main technical tool in the spectral analysis of sparse trees is the following theorem:
\begin{theorem}[Theorem 2.4 in \cite{breuer-CMP}] \label{decomposition}
Let $\Gamma$ be a normal rooted SH tree of type $\{L_n,k_n\}_{n=1}^\infty$.
Let 
\beq \label{M-n}
M_n= \left\{ \begin{array}{ll}
\prod_{j=1}^n k_j-\prod_{j=1}^{n-1}k_j & n>1\\
k_1-1 & n=1\\
1 & n=0.
\end{array} \right.
\eeq
Furthermore, let $R_0=0$ and $R_n=L_n+1,$ for $n\geq1$.
Then $\Delta$ is unitarily equivalent to a direct 
sum of Jacobi matrices, each operating on a copy of $\ell^2(\bbZ^+)$:
\beq \label{direct-sum-equivalence1}
\Delta \cong \oplus_{n=0}^\infty (\underbrace{J_n \oplus J_n \oplus \cdots \oplus J_n}_{M_n \textrm{ times}}) 
\eeq 
where $J_n=J(\{a_n(j)\}_{n=1}^\infty,\{b_n(j)\}_{n=1}^\infty)$ with 
\beq \label{a-n-j}
a_n(j)= \left\{ \begin{array}{ll}
\sqrt{k_m} & j=R_m-R_n \textrm{ for some } m>n\\
1 & \textrm{otherwise}
\end{array} \right.
\eeq
and
\beq \label{b-n-j}
b_n(j) \equiv 0.
\eeq
\end{theorem}
\begin{remarks}
\sloppy
1. The term - Jacobi matrix, with the notation $J(\{a(j)\}_{j=1}^\infty,\{b(j)\}_{j=1}^\infty)$, stands for the semi-infinite 
matrix 
\beq \label{jacobi}
J(\{a(j)\}_{j=1}^\infty,\{b(j)\}_{j=1}^\infty)=\left( \begin{array}{ccccc}
b(1)    & a(1) & 0      & 0      & \ldots \\
a(1)    & b(2) & a(2)    & 0      & \ldots \\
0      & a(2) & b(3)    & a(3)    & \ddots \\
\vdots & \ddots   & \ddots & \ddots & \ddots \\
\end{array} \right)
\eeq
with \beq \no  b(j) \in \bbR, \ a(j)>0. 
\eeq 
\fussy

2. For the case of a regular tree, a similar decomposition was discussed in \cite{froese, golenia, romanov} (see also 
\cite{solomyak} for a related result in the case of a metric tree).

3. As noted above, this theorem holds for $\ti{\Delta}$ as well, with the decomposition:
\beq \label{direct-sum-equivalence2}
\ti{\Delta} \cong \oplus_{n=0}^\infty (\underbrace{\ti{J}_n \oplus \ti{J}_n \oplus \cdots \oplus \ti{J}_n}_
{M_n \textrm{ times}})
\eeq
where $\ti{J_n}=J(\{\ti{a}_n(j)\}_{n=1}^\infty,\{\ti{b}_n(j)\}_{n=1}^\infty)$ with
\beq \label{tia-n-j}
\ti{a}_n(j)=a_n(j)
\eeq
and
\beq \label{tib-n-j}
\ti{b}_n(j)= \left\{ \begin{array}{ll}
-k_m-1 & j=R_m-R_n \textrm{ for some }m>n\\
-2 & \textrm{otherwise}
\end{array} \right.
\eeq
4. Note that each $J_n$ is a `tail' of $J_{n-1}$ in the sense that one can get $J_n$ by deleting a finite number of rows 
from the top, and the same number of columns from the left, of $J_{n-1}$.  
\end{remarks}

%%%%%%%%%%%%%%%%%%%%%%%%%%%%%%%%%%%%%%%%%%%%%%%%%%%%%%%%%%%%%%%%%%%% SECTION 3 %%%%%%%%%%%%%%%%%%%%%%%%%%%%%%%%%%%%%%

\section{Singular Measures on Sparse Trees}

In this and the next section we freely use terms (such as `transfer matrices') associated with the spectral theory 
of Jacobi matrices. The reader is referred to the appendix for their definitions, the notation we use, some basic results, 
and further references.

We start with a remark about the essential spectrum. Let $\Gamma$ be a sparse tree of type $\{L_n,k_n\}_{n=1}^\infty$. If 
$k_n \rightarrow \infty$, perturbation theory arguments show that the essential spectrum of $\Delta$ on $\Gamma$ is $[-2,2]$. If 
$\{k_n\}$ is bounded, then from remark 4 after theorem \ref{decomposition}, it is easy to see that the essential spectrum of $\Delta$ 
is contained in the essential spectrum of $J_0$. Since the reverse inclusion is immediate, we have the following

\begin{proposition} \label{essential-spectrum}
Let $\Gamma$ be a sparse tree of type $\{L_n,k_n\}_{n=1}^\infty$ and let $J_0=J_0(\Gamma)$ be the corresponding Jacobi
matrix appearing in theorem \ref{decomposition}. Let $\sigma_{\textrm{ess}}(\Delta)$ be the essential spectrum of 
$\Delta$
on $\Gamma$ and $\sigma_{\textrm{ess}}(J_0)$ be the essential spectrum of $J_0$. Then, if either $k_n \rightarrow \infty$
or $\{k_n\}$ is bounded, then
\beq \label{ess-spectrum-eq}
\sigma_{\textrm{ess}}(\Delta)=\sigma_{\textrm{ess}}(J_0).
\eeq
\end{proposition}

Now, let $H$ be a self-adjoint operator on a separable Hilbert space - $\mathcal{H}$, 
and $\psi \in \mathcal{H}$. The \emph{spectral measure} 
associated with $\psi$ and $H$ - $\mu_\psi$, is the unique measure on $\bbR$ satisfying
\beq \no
\langle \psi, (H-z)^{-1} \psi \rangle = \int_\bbR \frac{d\mu_\psi(x)}{x-z}, \qquad z \in \bbC \setminus \bbR
\eeq
(see e.g. \cite{reed-simon1}).

Theorem \ref{decomposition} reduces the spectral analysis of $\Delta$ on a sparse tree to the spectral 
analysis of Jacobi matrices with sparse `bumps' off the diagonal. (For $\ti{\Delta}$, we get `bumps' on the diagonal as 
well). An application of (suitably adapted) methods from the spectral theory of one-dimensional Schr\"odinger operators 
with sparse potentials (see \cite{last-review} for a review of the relevant theory) to 
the situation at hand, allows us to establish interesting 
spectral behavior for the Laplacian on certain sparse trees. 

The following basic lemma was proven in the appendix of \cite{aizenman} for the case of the Bethe lattice. It holds for 
any SH rooted tree.
\begin{lemma} \label{singular}
Let $\Gamma$ be a normal rooted SH tree with root - $O$. For any vertex $v$ in $\Gamma$, let $\delta_v \in \ell^2(\Gamma)$ be the 
delta function at $v$, and $\mu_v$ - the spectral measure associated with $\Delta$ and $\delta_v$. 
Let $(a,b)$ be an interval on which the absolutely continuous part of $\mu_O$ vanishes. Then, for any vertex $v$ of 
$\Gamma$, the absolutely continuous part of $\mu_v$ vanishes on $(a,b)$. 
\end{lemma}
\begin{remark}
Throughout this paper, `absolutely continuous' means absolutely continuous with respect to Lebesgue measure.
\end{remark}
\begin{proof}
The proof is a simple consequence of the identification of the essential support of the absolutely continuous spectrum 
with the set of energies for which the Green's function has positive imaginary part, combined with the recursion relation 
(see \cite{aizenman})
for the diagonal elements of the forward resolvents (these are the resolvents of $\Delta$ restricted to the various forward 
subtrees of $\Gamma$).
\end{proof}
For a normal rooted SH tree - $\Gamma$, let $\{J_n(\Gamma)\}_{n=0}^\infty$ be the Jacobi matrices appearing in the 
decomposition of $\Delta$ on $\Gamma$, given by theorem \ref{decomposition}. Let $\mu_n$ be the spectral measure associated 
with $\delta_1 \in \ell^2(\bbZ^+)$ and $J_n$. Then lemma \ref{singular} above, 
says that, if we want to prove that all spectral measures 
associated with the Laplacian on $\Gamma$ are singular, it suffices to prove that $\mu_0$ is singular. This will be useful 
later on. 

The following lemma is another simple tool for proving singularity of all the spectral measures. Its proof features the `bump' transfer 
matrix which will prove itself useful throughout the rest of this paper. (See equations \eqref{transfer-matrices1}-\eqref{t-j1-j2} for 
the definitions of the transfer matrices that we use below).
\begin{lemma} \label{unbounded-jacobi}
Let $R_m$ be a strictly increasing sequence of natural numbers such that, for $m$ large enough, $R_{m+1}-R_m \geq 2$.  
Let
\beq \no
J=J(\{a(j)\}_{j=1}^\infty,\{b(j)\}_{j=1}^\infty)
\eeq 
be a Jacobi matrix such that  
\beq \label{a-j}
a(j)= \left\{ \begin{array}{ll}
\rho_m & j=R_m \textrm{ for some } m\\
1 & \textrm{otherwise}
\end{array} \right.
\eeq
and
\beq \label{b-j}
b(j) \equiv 0,
\eeq
where $\rho_m>\delta>0$ for all $m$.
Then, if $\{\rho_m\}_{m=1}^\infty$ is unbounded, then the spectral measure, $\mu$, associated with 
$\delta_1 \in \ell^2(\bbZ^+)$ and $J$, is singular with respect to Lebesgue measure. 
\end{lemma}
\begin{proof}
Assume $\lim_{l\rightarrow\infty}\rho_{m_l}=\infty$, $R_{m_l}-R_{m_l-1}\geq 2$, $R_{m_l+1}-R_{m_l}\geq 2$ and fix 
$E \in \bbR$. Then 
\begin{align} \label{S'}
& S'_{m_l} \equiv T_{R_{m_l}+1,R_{m_l}-1}(E) = S_{R_{m_l}+1}(E)S_{R_{m_l}}(E) \notag \\
& =
\left( \begin{array}{cc}
\left(\frac{E^2}{\rho_{m_l}}-\rho_{m_l} \right) & -\frac{E}{\rho_{m_l}} \\
\frac{E}{\rho_{m_l}} & \frac{-1}{\rho_{m_l}} \end{array} \right). 
\end{align}
It follows that 
\beq \label{S'-norm}
\max(1,\rho_{m_l}-E^2) \leq \parallel T_{R_{m_l}+1,R_{m_l}-1}(E)^{-1} \parallel,
\eeq
and so, applying proposition \ref{last-simon} with $m_j=R_{m_l}+1$ and $l_j=R_{m_l}-1$, (note that $a_{R_{m_l}-1}=1$), 
we see that $\mu$ is singular on 
$\bbR$.    
\end{proof}

\begin{corollary} \label{singular-tree}
Let $\Gamma$ be a sparse tree of type $\{L_n,k_n\}_{n=1}^\infty$, with $\{k_n\}_{n=1}^\infty$ unbounded. Then all the spectral 
measures for $\Delta$ on $\Gamma$ are singular with respect to Lebesgue measure. 
\end{corollary}
\begin{proof}
The statement follows from theorem \ref{decomposition}, lemma \ref{unbounded-jacobi} and the fact that for a general Jacobi 
matrix, the vector $\delta_1$ is a cyclic vector.
\end{proof}

On the other hand, a simple consequence of proposition \ref{simon-stolz} is the following 
\begin{lemma} \label{continuous-jacobi}
Let $R_m$ be a strictly increasing sequence of natural numbers.  
Let 
\beq \no
J=J(\{a(j)\}_{j=1}^\infty,\{b(j)\}_{j=1}^\infty)
\eeq 
be a Jacobi matrix such that  
\beq \label{a-j2}
a(j)= \left\{ \begin{array}{ll}
\rho_m & j=R_m \textrm{ for some } m\\
1 & \textrm{otherwise}
\end{array} \right.
\eeq
and
\beq \label{b-j2}
b(j) \equiv 0,
\eeq
where $\rho_m>1$ for all $m$.
Let $\{\beta_m\}_{m=1}^\infty$ be a sequence such that $\beta_m \geq \rho_m$ and 
$\lim_{m \rightarrow \infty}\beta_m=\infty$ and let 
$A_m=\prod_{l=1}^m \beta_l^2$. If for some $\varepsilon>0$, 
\beq \label{sc-condition}
\limsup_{m \rightarrow \infty} \frac{(R_{m+1}-R_m)}{A_m^{(1+\varepsilon)}}>0
\eeq 
then the spectral measure, $\mu$, associated with 
$\delta_1 \in \ell^2(\bbZ^+)$ and $J$, is continuous on $(-2,2)$. 
\end{lemma} 
\begin{proof}
Consider an arbitrary closed interval $I \subseteq (-2,2)$. We will show that $\mu(I \cap \cdot)$ is continuous. 
From this it will follow that $\mu((-2,2) \cap \cdot)$ is continuous. 
For $R_m+1 < j < R_{m+1}$ we have that 
\beq \no
S_j(E)=\left( \begin{array}{cc}
E & -1 \\
1 & 0 \end{array} \right) 
\eeq
so that $\det (S_j(E))=1$. Furthermore, define, as in the proof of lemma \ref{unbounded-jacobi}
\begin{align} \label{S'1}
& S'_{m}(E) \equiv T_{R_m+1,R_m-1}(E)=S_{R_m+1}(E)S_{R_m}(E) \notag \\
&=
\left( \begin{array}{cc}
\left(\frac{E^2}{\rho_m}-\rho_m \right) & -\frac{E}{\rho_m} \\
\frac{E}{\rho_m} & \frac{-1}{\rho_m} \end{array} \right). 
\end{align}
Then we also have $\det (S'_{m}(E))=1$. Thus, if $j_2 \neq R_m$ and $j_1 \neq R_n$ for any $m$, $n$,  we have that 
\beq \label{unimodularity1}
\det (T_{j_1,j_2}(E))=1.
\eeq
Note that if $R_m+1 \leq j_2 < j_1 < R_{m+1}$, then $T_{j_1,j_2}(E)$ is just the transfer matrix for the free Laplacian, 
so that there is a constant $C_I$, depending only on the interval $I$, such that 
$1 \leq \parallel T_{j_1,j_2}(E) \parallel< C_I$ for any such $j_1,\ j_2$ and $E \in I$. In addition, we have from \eqref{S'1}
\beq \label{S'-norm-upper-bound}
\parallel S'_{R_{m}(E)} \parallel \leq \rho_m+5.
\eeq

Thus, for $R_m < j < R_{m+1}$ we have
\beq \label{upper-bound-for-T}
\parallel T_j(E) \parallel \leq C^m \prod_{j=0}^m \rho_j
\eeq
for some constant $C$ depending on $I$. 

Let $m_l \rightarrow \infty$ be a sequence for which 
\beq \label{subsequence-condition}
\frac{R_{m_l+1}-R_{m_l}}{A_{m_l}^{(1+\varepsilon)}}>\delta
\eeq
for some $\delta>0$. Let $M$ be chosen so that for all $m>M$, $\beta_m^\varepsilon>C$. Now, for sufficiently large $m_l>N$,
we have, from \eqref{upper-bound-for-T},
\beq \label{divergence-of-tail}
\sum_{j=R_{m_l}+1}^{R_{m_l+1}-1} \parallel T_j(E) \parallel^{-2} \geq 
\frac{\delta}{2}A_{m_l}^\varepsilon C^{-2m_l} \geq \frac{\delta}{2}C^{-2N}.
\eeq
Thus, the tail of the sum in \eqref{no-ev-condition} does not converge to zero. Therefore the sum is divergent and $\mu$ 
has no eigenvalues in $I$. This proves the lemma.
\end{proof}
\begin{corollary} \label{continuous-tree}
Let $\Gamma$ be a rooted SH tree of type $\{L_n,k_n\}_{n=1}^\infty$. Assume that, for some $\varepsilon>0$,
\beq \label{no-ev-condition2}
\limsup_{n \rightarrow \infty} \frac{(L_{n+1}-L_n)}{A_n^{(1+\varepsilon)}}>0,
\eeq
where $A_n=\prod_{l=1}^n \beta_l$ for some sequence $\beta_n \geq k_n$, such that $\beta_n \rightarrow \infty$. Then any 
spectral measure for $\Delta$ on $\Gamma$, is continuous on $(-2,2)$.
\end{corollary} 
\begin{proof}
The statement follows from theorem \ref{decomposition}, lemma \ref{continuous-jacobi} and the fact that for a general Jacobi 
matrix, the vector $\delta_1$ is a cyclic vector.
\end{proof}
A simple consequence of corollaries \ref{singular-tree}, \ref{continuous-tree} and proposition \ref{essential-spectrum} is the following theorem from 
\cite{breuer-CMP}:
\sloppy
\begin{theorem}[Theorem 4.1 in \cite{breuer-CMP}] \label{k-unbounded}
Let $\{k_n\}_{n=1}^\infty$ be a sequence of natural numbers such that $k_n \rightarrow \infty$ as $n \rightarrow \infty$.
Let $A_n=\prod_{j=1}^n k_j$. 
Assume that $(L_{n+1}-L_n) \rightarrow \infty$ and let $\Gamma$ be a SH rooted tree of type $\{L_n,k_n\}_{n=1}^\infty$.
Then the spectrum of $\Delta$ on $\Gamma$ consists of the interval $[-2,2]$ along with some discrete point spectrum 
outside this interval. If for some $\varepsilon>0$, 
\beq \label{sc-condition1}
\limsup_{n \rightarrow \infty} \frac{(L_{n+1}-L_n)}{A_n^{(1+\varepsilon)}}>0,
\eeq 
then any spectral measure for $\Delta$ is purely singular continuous on $(-2,2)$. 
\end{theorem}
\fussy

Since the next section discusses trees with bounded $k_n$, we quote the corresponding result from \cite{breuer-CMP}. We 
sketch its proof here since some of the ideas will appear in the sequel:
\begin{theorem}[Theorem 2.2 in \cite{breuer-CMP}] \label{k-bounded}
Let $k_0 \geq 2$ be a natural number and let $k_n \equiv k_0$. Assume that $(L_{n+1}-L_n) \rightarrow \infty$ and let 
$\Gamma$ be a SH rooted tree, of type $\{L_n, k_n\}_{n=1}^\infty$. 
Then the essential spectrum of $\Delta$ on $\Gamma$ contains the interval $[-2,2]$ and, provided 
$(L_{n+1}-L_n)$ increase sufficiently rapidly, any spectral measure for 
$\Delta$ is purely singular continuous on 
$(-2,2)$. By `sufficiently rapidly' we mean that $(L_{n+1}-L_n)$ has to be made sufficiently large with respect to 
$\{(L_{i+1}-L_i)\}_{i<n}$.
\end{theorem}

\begin{proof}
The claim about the essential spectrum follows immediately from proposition \ref{essential-spectrum}.

We want to show that if $(L_{n+1}-L_n)$ grow sufficiently rapidly (in the sense described in the theorem), then 
$\mu_{\delta_O}$ 
(the spectral measure of the delta function at the root of $\Gamma$) is purely singular continuous on $(-2,2)$. Lemmas 
\ref{singular} and \ref{continuous-jacobi} say that this suffices to imply that $(L_{n+1}-L_n)$ may be made to 
grow so fast as to make \emph{all} spectral measures singular continuous.  

Thus, our problem is reduced to the problem of studying a Jacobi matrix of the form $\Jab$ with 
\beq \label{a-j-sparse}
a(j)= \left\{ \begin{array}{ll}
\sqrt{k} & j=L_n+1 \textrm{ for some }n\\
1 & \textrm{otherwise}
\end{array} \right.
\eeq
and
\beq \label{b-j-sparse}
b(j) \equiv 0.
\eeq
The proof now follows closely Pearson's classical proof \cite{pearson}. 
For any $E \in (-2,2)$, let $\phi \in (0,\pi)$ be defined by $2\cos(\phi)=E$, and let  
$u_{1,2\cos(\phi)} \equiv u_{1,E}$ as defined in the appendix. 
Define the EFGP variables \cite{efgp} $r_\phi(j)$ and $\theta_\phi(j)$ through:
\beq \label{efgp1}
r_\phi(j)\cos(\theta_\phi(j))=u_{1,E}(j)-\cos(\phi)u_{1,E}(j-1)
\eeq
\beq \label{efgp2}
r_\phi(j)\sin(\theta_\phi(j))=\sin(\phi)u_{1,E}(j-1).
\eeq
It is easy to see that for 
$L_n+3 < j \leq L_{n+1}+1$, $r_\phi(j)=r_\phi(j-1)$ and $\theta_\phi(j)=\theta_\phi(j-1)+\phi$ (since the evolution equations for $u_{1,E}$ 
there coincide with those of the usual Laplacian on $\bbZ^+$), so that for all $n$, 
\beq \label{constant-efgp1}
r_\phi(L_n+1)=r_\phi(L_{n-1}+3)
\eeq 
and
\beq \label{constant-efgp2} 
\theta_\phi(L_n+1)=\theta_\phi(L_{n-1}+3)+(L_n-L_{n-1}-2)\phi.
\eeq

Since the map $g(\phi)=2\cos{\phi}$ is continuously invertible on $(0,\pi)$, 
it is clear that, instead of studying $\mu_{\delta_O}$ on a given closed interval in $(-2,2)$, 
we may study its push-forward (via $g^{-1}$) - $\nu$ - on the corresponding closed interval in $(0,\pi)$. 
Now, since $a(j)$ is bounded, standard methods (see e.g.\ \cite{carmona} section III.3) show that $\nu$ on $(0,\pi)$ is equivalent to 
the measure defined by the limit
\beq \label{measure-asymptotics}
\lim_{n \rightarrow \infty} \int_I \frac{d\phi}{\left( r_\phi(L_n+3) \right)^2}
\eeq
where the integral is performed over closed subintervals $I \subseteq (0,\pi)$. Noting that 
\beq \no
u_{1,2\cos(\phi)}(j-1)=\frac{\sin(\theta_\phi(j))}{\sin(\phi)}r_\phi(j),
\eeq
\beq \no
u_{1,2\cos(\phi)}(j)=\frac{\sin(\phi+\theta_\phi(j))}{\sin(\phi)}r_\phi(j),
\eeq
one can employ the transfer matrix 
\beq \no
S'_n(2\cos(\phi)) \equiv S_{L_n+2}(2\cos(\phi))S_{L_n+1}(2\cos(\phi))
\eeq 
and \eqref{constant-efgp1}-
\eqref{constant-efgp2} to express $(r_\phi(L_n+3))^2$ as a function of $(r_\phi(L_{n-1}+3))^2$. Note that $S'_n(2\cos(\phi))$ is 
unimodular. We get that 
\beq \label{pearson1}
\int_I \frac{d\phi}{r_\phi(L_n+3)^2}=\int_I d\phi\prod_{i=1}^n f_i(\phi,N_i\phi,\theta_{i-1}(\phi)),
\eeq
where $N_i=L_i-L_{i-1}-2$ and 
\beq \label{pearson2}
f_i(\phi,y,\theta)^{-1}=A(\phi)+B(\phi)\cos2(\theta+y)+C(\phi)\sin2(\theta+y),
\eeq
with $A^2-B^2-C^2=\det(S'_i(E))=1$. 

The situation described in \eqref{pearson1}-\eqref{pearson2} is exactly the same as that in section 3 of \cite{pearson} 
(although the explicit expressions one gets for $A$, $B$ and $C$ are different). Note that the precise form of $S'_n(E)$ is 
of no 
importance. The unimodularity of this transfer matrix, together with the fact that $a(L_n+1)=\textrm{const.} \neq 1$, 
suffice to imply that (as in \cite{pearson}) the corollary to theorem 1 from \cite{pearson}, applies in this situation. 
This shows that if the $L_n$ are chosen to increase rapidly enough, $\nu$, and thus $\mu_{\delta_O}$, is singular continuous.  
\end{proof}

%%%%%%%%%%%%%%%%%%%%%%%%%%%%%%%%%%%%%%%%%% SECTION 4 %%%%%%%%%%%%%%%%%%%%%%%%%%%%%%%%%%%%%%%%%%%%%%

\section{Finite Dimensional Trees}

As noted in the introduction, aside from providing interesting examples for the Laplacian, sparse trees are interesting as objects interpolating between 
the one dimensional line ($\bbZ^+$) and the Bethe lattice (which is infinite dimensional in a natural sense). By 
tuning the sequences $\{L_n\}$ and $\{k_n\}$, one may construct trees with dimensions having any real value between one 
and infinity. A simple example is obtained as follows:  
Let $k_n\equiv k \geq 2$ and take $L_n=[\gamma^n]$ for some $\gamma>1$. Denote the SH rooted tree of type 
$\{L_n,k_n\}_{n=1}^\infty$ by $\Gamma_{k,\gamma}$. A simple calculation gives:
\begin{proposition} \label{dimension}
Fix $\gamma>1$ and $\bbN \ni k \geq 2$. Let $\Gamma=\Gamma_{k,\gamma}$ and let 
$S_\Gamma(r) = \{v_i \in \mathscr{V}(\Gamma) \mid d(v_i,O) \leq r \}$, where $\mathscr{V}(\Gamma)$ is the set of vertices of $\Gamma$.
Then
\beq \no
\limsup_{r \rightarrow \infty} \frac{\log \#S_\Gamma(r)}{\log r}=
\liminf_{r \rightarrow \infty} \frac{\log \#S_\Gamma(r)}{\log r}=\frac{\log \gamma k}{\log \gamma}.
\eeq
\end{proposition}
Below, we shall refer to the quantity $\frac{\log \gamma k}{\log \gamma}$ as the \emph{dimension} of $\Gamma_{k,\gamma}$.

In the context of the analogy between sparse trees and Schr\"odinger operators with sparse potentials, described in the 
previous section, $\Gamma_{k,\gamma}$ is analogous to a Schr\"odinger operator with bumps of fixed height placed at the 
sites $[\gamma^n]$ of $\bbZ^+$. Zlato\v s deals with such operators in \cite{zlatos} and the analysis we present below is 
an adaptation 
of his methods (in particular - section 6 of \cite{zlatos}) to the case at hand. While a large part of the argument 
translates word for word, there are a few significant changes, mainly having to do with the fact that the transfer matrices 
for our case are not, in general, unimodular. This is important for some of the arguments and, therefore, has to be bypassed 
to get the same results here. We discuss the changes below and give a sketch of the proof. However, we refer the reader to 
\cite{zlatos} for a more detailed discussion.

First,
\begin{definition}
Let $\Gamma$ be a rooted tree. For any self-adjoint operator $H$ on
$\ell^2(\Gamma)$,
and $-\frac{\pi}{2}<\varrho<\frac{\pi}{2}$, let
\beq \label{H-varrho}
H_\varrho=H-\tan(\varrho)P_O
\eeq
where $P_O$ is the orthogonal projection onto the subspace spanned by the delta
function at $O$.
We refer to $H_\varrho$ as $H$ with boundary condition - $\varrho$.
\end{definition}

We also need:
\begin{definition}
Let $\mu$ be a measure on $\bbR$. We say that $\mu$ has \emph{exact local dimension} in $I \subseteq \bbR$ if for any 
$E \in I$ there is an $\alpha(E)$ and for any $\varepsilon>0$ there is $\delta>0$ for which 
$\mu((E-\delta,E+\delta)\cap \cdot)$ is both continuous with respect to $(\alpha(E)-\varepsilon)$-dimensional 
Hausdorff measure, and singular with respect to $(\alpha(E)+\varepsilon)$-dimensional Hausdorff measure. We call 
$\alpha(E)$ the \emph{local dimension} of the measure $\mu$.
\end{definition}
Let $\omega_n$ be a random variable
uniformly distributed over 
\beq \no
[-n,-n+1,\ldots,n-1,n].
\eeq 
Let $(\Omega,\bbP)$ be the product
probability space for all $\omega_n$, $n=1,2,\ldots$. Fix $1<k\in \bbN$ and $\gamma>2$ and for each
$\omega \in \Omega$ let $\Gamma_{k,\gamma}^\omega$ be the SH rooted tree of type $\{L_n^\omega,k_n\}_{n=1}^\infty$ 
for $L_n^\omega=[\gamma^n]+\omega_n$ and $k_n \equiv k$. Clearly, proposition \ref{dimension} holds for any 
$\Gamma_{k,\gamma}^\omega$. The main result of this section is

\begin{theorem} \label{finite-dimensional-trees}
For $\bbP$-a.e.\ $\omega$, all the spectral measures for $\Delta$ on $\Gamma_{k,\gamma}^\omega$ are singular with respect 
to Lebesgue measure. 
Furthermore, let $V(k)=\frac{(1+k)^2}{4k}$ and let
\beq \label{I-definition}
I=\left( -\sqrt{\frac{8(\gamma-V(k))}{2\gamma-1}},\sqrt{\frac{8(\gamma-V(k))}{2\gamma-1}} \right)
\eeq
if $\gamma>V$, and $I=\emptyset$ otherwise. Then
for $\bbP$-a.e.\ $\omega$ and for Lebesgue a.e.\ $\varrho \in (-\frac{\pi}{2},\frac{\pi}{2})$, 
the spectral measure, 
$\mu_{\delta_O}$, 
associated with $\Delta_{\varrho}$ on $\Gamma_{k,\gamma}^\omega$, and with the delta function at the root, is purely 
singular continuous in $I$ with exact local dimension
\beq \label{exact-Hausdorff-dimension}
1-\frac{\log(\frac{4V(k)-\frac{E^2}{2}}{4-E^2})}{\log(\gamma)}
\eeq
and it is dense pure point in the rest of $[-2,2]$.
\end{theorem}

\begin{corollary}
Assume $\gamma \geq 4$. Then if the dimension of $\Gamma_{k,\gamma}^\omega$ is at least $3$, we have that $I=\emptyset$ and 
so, for $\bbP$-a.e.\ $\omega$ and for Lebesgue a.e.\ $\varrho \in (-\frac{\pi}{2},\frac{\pi}{2})$, the spectral measure 
$\mu_{\delta_O}$ is dense pure point in $[-2,2]$.
\end{corollary}

\begin{proof}[Proof of the corollary]
This is a simple computation.
\end{proof}

\begin{proof}[Proof of theorem \ref{finite-dimensional-trees}]
Theorem \ref{decomposition} and lemma \ref{singular} imply that the theorem is an immediate consequence of proposition 
\ref{random-an} below. 
\end{proof}

\begin{proposition} \label{random-an}
Fix $k \geq 2$ and $\gamma>2$.
For any $\omega \in \Omega$ and $-\frac{\pi}{2}<\varrho<\frac{\pi}{2}$, let 
$J^\omega_\varrho=J(\{a^\omega(j)\},\{b_\varrho(j)\})$ be a Jacobi matrix with 
\beq \label{a-j-bounded}
a^\omega(j)= \left\{ \begin{array}{ll}
\sqrt{k} & j=[\gamma^m]+\omega_m+1 \textrm{ for some } m\\
1 & \textrm{otherwise}
\end{array} \right.
\eeq
and
\beq \label{b-j-bounded}
b_\varrho(j)= \left\{ \begin{array}{ll}
-\tan(\varrho) & j=1 \\
0 & \textrm{otherwise}
\end{array} \right. 
\eeq
Then, for $\bbP$-a.e.\ $\omega$, and for any $\varrho \in (-\frac{\pi}{2},\frac{\pi}{2})$, the spectral measure $\mu_\varrho$, 
associated 
with the vector $\delta_1 \in \bbZ^+$ and with the Jacobi matrix $J^\omega_\varrho$, is singular with respect to Lebesgue 
measure.
Furthermore, for $\bbP$-a.e.\ $\omega$ and for Lebesgue a.e.\ $\varrho \in (-\frac{\pi}{2},\frac{\pi}{2})$, the spectral 
measure 
$\mu_\varrho$ is purely singular continuous in $I$ with exact local dimension given by \eqref{exact-Hausdorff-dimension}, 
where $I$ is as defined in theorem \ref{finite-dimensional-trees}, and it is dense pure point in the rest of $[-2,2]$. 
\end{proposition}

A central role in the proof of the proposition will be played by the EFGP transform introduced in the proof of theorem 
\ref{k-bounded}: 

Fix $k \geq 2$ and let $J=\Jab$ be a Jacobi matrix satisfying 
\beq \label{a-j-sparse1}
a(j)= \left\{ \begin{array}{ll}
\sqrt{k} & j=L_n+1 \textrm{ for some }n\\
1 & \textrm{otherwise}
\end{array} \right.
\eeq
for a sequence $\{L_n\}$ satisfying $L_{n+1}-L_n \geq 2$,
and
\beq \label{b-j-sparse1}
b(j) \equiv 0.
\eeq 
For any $E \in (-2,2)$, let $\phi \in (0,\pi)$ be defined by $2\cos(\phi)=E$, and let  
$u_{2\cos(\phi)} \equiv u_{E}$ solve \eqref{schrodinger3} for $J$, namely
\beq \label{schrodinger4}
a(j)u(j+1)+a(j-1)u(j-1)=Eu(j), \qquad j \geq 1
\eeq 
with $a(0)=1$. 
Recall that the EFGP variables \cite{efgp} corresponding to $u$, $r_\phi(j)$ and $\theta_\phi(j)$, are defined through:
\beq \label{mefgp1}
r_\phi(j)\cos(\theta_\phi(j))=u_{E}(j)-\cos(\phi)u_{E}(j-1)
\eeq
\beq \label{mefgp2}
r_\phi(j)\sin(\theta_\phi(j))=\sin(\phi)u_{E}(j-1).
\eeq
First, note that there are positive constants, $C_1(\phi)$, $C_2(\phi)$, such that
\beq \label{rj-bound-for-u}
C_1(\phi)(|u(j-1)|^2+|u(j)|^2)\leq r_\phi(j)^2 \leq C_2(\phi)(|u(j-1)|^2+|u(j)|^2).
\eeq
We call $r(j)$ the \emph{EFGP norm} of $u(j)$. Note, also, that for any $j$
\beq \label{rj-bound-for-u2}
C_1(\phi)(\sum_{i=(j-1)}^{(j+2)}|u(i)|^2)\leq r_\phi(j)^2+r_\phi(j+2)^2 \leq C_2(\phi)(\sum_{i=(j-1)}^{(j+2)}|u(i)|^2).
\eeq
For a function $f: \bbZ^+ \rightarrow \bbC$ and a sequence $\mathcal{L}=\{L_n\}_{n=1}^\infty$ of natural numbers, define
\beq \no
\parallel f \parallel_{L,\mathcal{L}}=
\left( \sum_{1 \leq j \leq L, j \neq L_n+2 \textrm{ for any }n} |f(j)|^2 \right)^{1/2}
\eeq
for any $L \in \bbN$, and extend to $L \in \bbR$ by linear interpolation.
Then we have 

\begin{lemma} \label{m-jit-last2}
Let $J=\Jab$ be as defined in \eqref{a-j-sparse1}-\eqref{b-j-sparse1} and let $\mathcal{L}=\{L_n\}$. 
Assume that for some $1 \leq \beta <2$ and every $E$ in some Borel set $A \subseteq (-2,2)$, the EFGP norm, $r_u$, of 
every solution $u$ of \eqref{schrodinger4} obeys
\beq \label{m-jit-last-condition2}
\limsup_{L \rightarrow \infty} \frac{\parallel r_u \parallel_{L,\mathcal{L}}^2}{L^\beta}<\infty.
\eeq  
Then $\mu(A \cap \cdot)$ is continuous with respect to $(2-\beta)$-dimensional Hausdorff measure, where $\mu$ is the 
spectral measure associated with $\delta_1$ and $J$.
\end{lemma}

\begin{lemma} \label{m-jit-last3}
Let $J=\Jab$ be as defined in \eqref{a-j-sparse1}-\eqref{b-j-sparse1} and let $\mathcal{L}=\{L_n\}$. Let $u_{1,E}$ be the solution of 
\eqref{schrodinger4} that satisfies 
\beq \no
u(0)=0,\ u(1)=1,
\eeq 
and let $r_{1,E}$ denote the corresponding EFGP norm. If  
\beq \label{m-jit-last-condition3}
\liminf_{L \rightarrow \infty} \frac{\parallel r_{1,E} \parallel_{L,\mathcal{L}}^2}{L^\alpha}=0
\eeq
for every $E$ in some Borel set $A \subseteq (-2,2)$, then $\mu(A \cap \cdot)$ is singular with respect to 
$\alpha$-dimensional Hausdorff measure, where $\mu$ is the spectral measure associated with $\delta_1$ and $J$.
\end{lemma}

\sloppy
\begin{proof}[Proof of lemmas \ref{m-jit-last2} and \ref{m-jit-last3}]
The proofs follow immediately from \eqref{rj-bound-for-u2}, from the fact that $L_{n+1}-L_n \geq 2$, 
and from propositions \ref{jit-last2} and 
\ref{jit-last3} respectively.
\end{proof}
\fussy

Our aim, therefore, is to control the growth of the EFGP norms of generalized eigenfunctions.
We recall from the previous section that for  
$L_n+3 < j \leq L_{n+1}+1$, $r_\phi(j)=r_\phi(j-1)$ and $\theta_\phi(j)=\theta_\phi(j-1)+\phi$ and that, therefore,  
\beq \label{mconstant-efgp1}
r_\phi(L_n+1)=r_\phi(L_{n-1}+3)
\eeq 
and
\beq \label{mconstant-efgp2} 
\theta_\phi(L_n+1)=\theta_\phi(L_{n-1}+3)+(L_n-L_{n-1}-2)\phi.
\eeq
From lemmas \ref{m-jit-last2} and \ref{m-jit-last3}, we see that $r_\phi(L_n+2)$ is irrelevant so all we need is 
to find the relation between $r_\phi(L_n+1)$ and $r_\phi(L_n+3)$. 
As before, using the relations
\beq \no
u_{2\cos(\phi)}(j-1)=\frac{\sin(\theta_\phi(j))}{\sin(\phi)}r_\phi(j),
\eeq
\beq \no
u_{2\cos(\phi)}(j)=\frac{\sin(\phi+\theta_\phi(j))}{\sin(\phi)}r_\phi(j),
\eeq
and the transfer matrix  
\begin{align} \no
& S'_n(2\cos(\phi)) \equiv S_{L_n+2}(2\cos(\phi))S_{L_n+1}(2\cos(\phi)) \notag \\
& =\left( \begin{array}{cc}
\left(\frac{4\cos^2(\phi)}{\sqrt{k}}-\sqrt{k} \right) & -\frac{2\cos(\phi)}{\sqrt{k}} \\
\frac{2\cos(\phi)}{\sqrt{k}} & \frac{-1}{\sqrt{k}} \end{array} \right), \notag 
\end{align}
we get that 
\begin{align} \label{r-recursion}
& r_\phi(L_n+3)^2=r_\phi(L_n+1)^2 \cdot \notag \\ 
& \cdot \Big( A(\phi)+B(\phi)\cos(2\theta_\phi(L_n+1))+C(\phi)\sin(2\theta_\phi(L_n+1)) \Big),
\end{align}
where 
\beq \label{Aphi}
A(\phi)=\frac{1}{\sin^2(\phi)} \left(\frac{1+k^2}{2k}-\cos^2(\phi) \right) \geq 0
\eeq
and 
\beq \label{unimodul}
A(\phi)^2-B(\phi)^2-C(\phi)^2=1.
\eeq

Explicit formulas for $B(\phi)$ and $C(\phi)$ can be derived but they are of no consequence. The last relation follows 
from the fact that 
\beq \no
\det (S'_n(2\cos(\phi)))=1.
\eeq

\begin{proof} [Proof of proposition \ref{random-an}]
For a fixed $\phi \in (0,\pi)$ (and $E=2\cos(\phi) \in (-2,2)$), let 
\beq \label{f}
f(\theta)=\frac{1}{2}\log \Big( A(\phi)+B(\phi)\cos(2\theta)+C(\phi)\sin(2\theta) \Big). 
\eeq  
Now, for a given $\omega \in \Omega$, let $\mathcal{L}^\omega=\{L_n^\omega\}_{n=1}^\infty$ (recall  
$L_n^\omega=[\gamma^n]+\omega_n$). Let $u^\omega_{1,E}$ be a generalized 
eigenfunction for $J^\omega_0$ such that $u^\omega_{1,E}(0)=0$ and 
$u^\omega_{1,E}(1)=1$, and let $r^\omega_{1,\phi}$ be the corresponding EFGP norm.  
Then, the above implies that 
\begin{align} \label{rLn+3}
& Y_n(\omega) \equiv  \log r^\omega_{1,\phi}(L_n^\omega+3)-\log r^\omega_{1,\phi}(L_{n-1}^\omega+3) \notag \\
& =\frac{1}{2} \log \Big( A(\phi)+B(\phi)\cos(2\theta_\phi(L_n^\omega+1))+C(\phi)\sin(2\theta_\phi(L_n^\omega+1)) \Big). 
\end{align}
Note that, if $\omega,\eta \in \Omega$ are such that $\omega_j=\eta_j$ for $j=1,2,\ldots,n-1$, and 
$\eta_n=\omega_n+l$ for some $l$, and $Y_n(\omega)=f(\theta)$ for some $\theta$, then 
$Y_n(\eta)=f(\theta+l \phi)$.

Let
\begin{align} \label{Zomega}
& Z=\frac{1}{2\pi}\int_0^{2\pi}f(\theta)d\theta=\frac{1}{2}\log \frac{A(\phi)+\sqrt{A(\phi)^2-B(\phi)^2-C(\phi)^2}}{2} \notag \\
& =\frac{1}{2}\log \frac{A(\phi)+1}{2}=\frac{1}{2}\log \left( \frac{V(k)}{ \sin^2(\phi)}-\frac{\cot^2(\phi)}{2}\right),
\end{align}
and 
\beq \label{ftilde}
\ti{f}(\theta)=f(\theta)-Z.
\eeq
Then $\ti{f}(\theta)=\ti{f}(\theta+2\pi)$, $\int_0^{2\pi}\ti{f}(\theta)d\theta=0$ and $\ti{f}$ is $C^4$, so $\ti{f}$ 
satisfies the conditions of lemma 6.2 in \cite{zlatos}. Thus, it follows that, for any $\phi$ in a set of full Lebesgue measure - 
$S \subseteq (0,\pi)$, and for $\bbP$-a.e.\ $\omega$,
\beq \no
\frac{\sum_{n=1}^N (Y_n(\omega)-Z)}{N} \rightarrow 0,
\eeq
(see section 6 of \cite{zlatos}). From this we get that for any $\phi \in S$ and for a.e.\ $\omega \in \Omega$, the EFGP 
norm $r^\omega_{1,\phi}$ satisfies
\beq \label{bounds-for-romega}
(L_n^\omega+3)^{d_1} \leq r^\omega_{1,\phi}(L_n^\omega+3) \leq (L_n^\omega+3)^{d_2} 
\eeq 
for any 
\beq \no
d_1 < \frac{Z}{\log \gamma} < d_2
\eeq
provided $n$ is large enough. By Fubini's theorem, we get that \eqref{bounds-for-romega} holds for a.e.\ $\omega$ and Lebesgue almost 
every $\phi$.
Furthermore, lemma \ref{subordinate-existence} below assures us that for each such $(\omega,\phi)$, there exists a subordinate 
solution $u^\omega_{\textrm{sub},\phi}$ to \eqref{schrodinger4} whose EFGP norm $r^\omega_{\textrm{sub},\phi}$ satisfies
\beq \label{subordinate-romega}
r^\omega_{\textrm{sub},\phi}(L_n^\omega+3) \leq (L_n^\omega+3)^{-d_1}  
\eeq
for large enough $n$. Since the absolutely continuous part of $\mu_\varrho$ (for any $\varrho$) is supported on the set of 
energies with no 
subordinate solution, it follows that, for any $\varrho$, $\mu_\varrho$ is purely singular on $(-2,2)$. Furthermore, since 
the singular part of $\mu_\varrho$ is supported on the set of energies where the subordinate solutions obey the 
appropriate boundary conditions, rank-one perturbation arguments (see e.g.\ \cite{rankone}), imply that for Lebesgue 
almost any $\varrho$, $\mu_\varrho$ is supported on the set of energies where the subordinate solutions behave as in 
\eqref{subordinate-romega}. 
Noting that, in this case, 
\beq \no
\parallel r^\omega_{\textrm{sub},\phi} \parallel_{L,\mathcal{L}^\omega}^2 \leq C L^{1-2\frac{Z}{\log \gamma}}
\eeq
if $\frac{Z}{\log \gamma} \leq \frac{1}{2}$($\Leftrightarrow E \in I$) and is square-summable otherwise, and that 
the EFGP norm - $r^\omega_\phi$ of any other solution to $\eqref{schrodinger4}$ satisfies 
\beq \no
\parallel r^\omega_\phi \parallel_{L,\mathcal{L}^\omega}^2 \leq C L^{1+2\frac{Z}{\log \gamma}},
\eeq
we finish the proof with the help of lemmas \ref{m-jit-last2} and \ref{m-jit-last3}.
\end{proof}

The following lemma was used in the proof above. It is essentially the same as lemma 2.1 from \cite{zlatos}, the only 
difference being in that not all transfer matrices for the general Jacobi case are unimodular. Since not all transfer 
matrices are relevant for the argument, this is inconsequential. In the proof below we only demonstrate this point. 
For the complete proof the reader is referred to \cite{zlatos}. 
\begin{lemma} \label{subordinate-existence}
Let $J=\Jab$ be a Jacobi matrix satisfying 
\beq \label{a-j-sparse2}
a(j)= \left\{ \begin{array}{ll}
\sqrt{k} & j=L_n+1 \textrm{ for some }n\\
1 & \textrm{otherwise}
\end{array} \right.
\eeq
for a sequence $\{L_n\}$ satisfying $L_{n+1}-L_n \geq 2$,
and
\beq \label{b-j-sparse2}
b(j) \equiv 0.
\eeq
Assume that for some $E \in (-2,2)$, $u$ is a solution of \eqref{schrodinger4} whose EFGP norm satisfies
\beq \label{EFGP-norm-growth}
r(L_n+3)=e^{\sigma_n}
\eeq
where $\sigma_n=\sum_{j=1}^n (Z_j+X_j)$ with $0<d_1 \leq Z_j \leq d_2 < \infty$ and $\sum_{j=1}^n X_j=o(n)$.
Then there exists a subordinate solution $v$ of \eqref{schrodinger4} for $E$, such that for any $d<d_1$ and for 
all sufficiently large $n$, the corresponding EFGP norm - $p$ - satisfies
\beq \label{EFGP-norm-decay}
p(L_n+3) \leq e^{-dn}.
\eeq
\end{lemma}
\begin{proof}
Let $v$ be any solution of \eqref{schrodinger4} different from $u$ and let $p$ be its EFGP norm. Since 
the transfer matrices $T_{L_n+2}(E)$ are unimodular, the argument of theorem 2.3 from \cite{efgp} applies to 
show that there exist $E$-dependent constants, $c_1,c_2$ such that 
\beq \label{transfer-matrix-EFGP}
c_1 \max (p_{L_n+3},r_{L_n+3}) \leq \parallel T_{L_n+2}(E) \parallel \leq c_2 \max (p_{L_n+3},r_{L_n+3}). 
\eeq
Note, further, that there exists a constant $B>0$ such that 
\beq \label{transfer-matrix-bound}
\parallel T_{n+1,n}(E) \parallel \equiv \parallel T_{L_{n+1}+2,L_n+2}(E) \parallel \leq B
\eeq
and that $\det (T_{n+1,n}(E))=1$ as well. Thus, it follows that 
\beq \label{simon-last-condition}
\sum_{n=1}^\infty \frac{\parallel T_{n,n-1}(E) \parallel^2}{\parallel T_{L_n+2}(E) \parallel^2}<\infty
\eeq
so one can apply theorem 8.1 from \cite{last-simon} to get a vector $\overline{v} \in \mathbb{R}^2$ such that
\beq \no
\frac{\parallel T_{L_n+2}(E) \overline{v} \parallel}{\parallel T_{L_n+2}(E) \overline{u} \parallel} \rightarrow 0
\eeq
for any other vector $\overline{u} \in \bbR^2$. From this point, the proof follows the proof of lemma 2.1 from 
\cite{zlatos}, word for word, to show that $\overline{v}$ generates the claimed solution.
\end{proof}

%%%%%%%%%%%%%%%%%%%%%%%%%%%%%%%% APPENDIX %%%%%%%%%%%%%%%%%%%%%%%%%%%%%%%%%%%%%%%%%%
\appendix
\section{Self-Adjointness of the Laplacian on Normal SH Rooted Trees}

\begin{proposition} \label{self-adjointness}
Let $\Gamma$ be a rooted SH tree of type $\{L_n,k_n\}_{n=1}^\infty$. Then the operator $\Delta_1$ defined over 
\beq \label{domain-of-delta1}
D(\Delta_1)=\{u \in \ell^2(\Gamma) \mid u \textrm{ is of compact support} \}
\eeq
via the equation
\beq \label{formal-lap1}
(\Delta_1 u)(x)=\sum_{y : d(x,y)=1}u(y),
\eeq
is symmetric. If $\limsup_{n \rightarrow \infty} (L_{n+1}-L_n) >1$ then $\Delta$ - the closure of $\Delta_1$ - 
is self-adjoint. The same statement holds for $\ti{\Delta}_1$ and $\ti{\Delta}$ (defined over the same domain), 
with equation \eqref{formal-lap1} replaced by
\beq \label{formal-lap2}
(\ti{\Delta}_1 u)(x)=\sum_{y : d(x,y)=1}u(y)-\#\{y : d(x,y)=1 \}\cdot u(x).
\eeq
\end{proposition} 

\begin{proof}
Since the proof for $\Delta$ and $\ti{\Delta}$ is precisely the same, we use $\Delta$. 

It is trivial to see that $\Delta_1$ is symmetric, so in order to show that $\Delta$ is self-adjoint, all we have to show is 
that $\ker(\Delta^* \pm i)=\{0\}$.

Assume that $\Delta u=i u$, then it follows that $\Delta \overline{u}= -i \overline{u}$ (where 
$\overline{(\cdot)}$ for a complex number denotes 
complex conjugation). Let $n_j$ be a subsequence for which $L_{n_j}+1<L_{n_j+1}$. For a vertex $v$ with 
$|v| \equiv d(v,O)=L_{n_j}+1$, 
let us denote its \emph{unique} forward neighbor by $\hat{v}$. One can verify that an analogue of Green's 
formula (see e.g.\ \cite{berez} Chapter VII, formula 1.4) holds and we have:
\beq \no 
2i\sum_{|v|\leq L_{n_j}+1} |u(v)|^2=\sum_{|v|=L_{n_j}+1}u(\hat{v}) \cdot 
\overline{u(v)}-\sum_{|v|=L_{n_j}+1} \overline{u(\hat{v})} \cdot u(v)
\eeq  
so
\beq \no
\sum_{|v|\leq L_{n_j}+1} |u(v)|^2 \leq \sum_{|v|=L_{n_j}+1}|u(\hat{v})| \cdot
|u(v)|. 
\eeq
From this it follows that if $u \neq 0$, then the RHS above does not converge to zero and therefore 
$u \notin \ell^2(\Gamma)$. This proves the proposition.
\end{proof}

%%%%%%%%%%%%%%%%%%%%%%%%%%%%%%%%%%%%%%% APPENDIX 2 %%%%%%%%%%%%%%%%%%%%%%%%%%%%%%%%%%%%%%%%%%%%%
\section{Eigenfunctions and Transfer Matrices for Jacobi Matrices}

Let $J=\Jab$ be a Jacobi matrix with $b(j) \in \bbR$ and $a(j)>0$ satisfying $\sum_{j=1}^\infty \frac{1}{a(j)}= \infty$ 
(which suffices for $J$ to be self-adjoint \cite{berez}). A basic idea in the spectral theory of Jacobi matrices is to 
relate spectral properties of $J$ as reflected by $\mu=\mu_{\delta_1}$ - the spectral measure of $\delta_1$ - to 
properties of formal eigenfunctions
\beq \label{formal-ev}
Ju=Eu.
\eeq
By this term we mean functions $u: \bbZ^+ \rightarrow \bbC$ which satisfy
\beq \label{schrodinger1}
a(j)u(j+1)+a(j-1)u(j-1)+b(j)u(j)=Eu(j), \qquad j>1
\eeq
\beq \label{schrodinger2}
a(1)u(2)+b(1)u(1)=Eu(1).
\eeq
Since, for a given 
$E \in \bbR$, all solutions to \eqref{schrodinger1}-\eqref{schrodinger2} are linearly dependent (determined by $u(1)$), 
it suffices to study 
$u_{1,E}$ which is the solution satisfying $u_{1,E}(1)=1$.
It is convenient to define $a(0)=1$ and to extend \eqref{schrodinger1} to $j=1$ by demanding $u_{1,E}(0)=0$. Thus $u_{1,E}$ 
is the unique solution to 
\beq \label{schrodinger3}
a(j)u(j+1)+a(j-1)u(j-1)+b(j)u(j)=Eu(j), \qquad j \geq 1
\eeq
with
\beq \label{u1}
u_{1,E}(0)=0, \qquad u_{1,E}(1)=1.
\eeq
We further define $u_{2,E}$ as the unique solution to \eqref{schrodinger3} satisfying 
\beq \label{u2}
u_{2,E}(0)=1, \qquad u_{2,E}(1)=0.
\eeq
Note that any solution, $u$, to \eqref{schrodinger3} with $u(1) \neq 0$ can be viewed as $u_{1,E}$ for a slightly modified 
Jacobi matrix. Namely, $u$ solves \eqref{schrodinger1}-\eqref{schrodinger2} for the same set of parameters except with 
$b(1)$ changed to $(b(1)+\frac{u(0)}{u(1)})$. This remark is basic for the analysis of section 4.

We say that $u_{1,E}$ is subordinate if 
\beq \label{subord1}
\lim_{L \rightarrow \infty} \frac{\parallel u_{1,E} \parallel_L}{\parallel u_{2,E} \parallel_L}=0
\eeq
where
\beq \no
\parallel f \parallel_L=\left( \sum_{j=1}^{[L]} |f(j)|^2+(L-[L])|f([L]+1)|^2 \right)^{1/2}.
\eeq
The Gilbert-Pearson theory of subordinacy \cite{subord} says that the singular part of $\mu$ is supported on the 
set of energies where 
$u_{1,E}$ is subordinate, and that the absolutely continuous part of $\mu$ is supported off this set. 
The Jitomirskaya-Last extension of this theory \cite{jit-last} analyzes further the 
singular part of $\mu$ according to its singularity/continuity with respect to dimensional Hausdorff measures (see \cite{jit-last} 
for the concept of Hausdorff measures and dimensions).  
In section 4 of the paper we use the following results from \cite{jit-last}:  
\begin{proposition}[\cite{jit-last}] \label{jit-last2}
Assume that $1 \leq a(j) \leq M$ for some $M>1$. Assume that for some $1 \leq \beta <2$ and every $E$ in some Borel set $A$, 
every solution $u$ of \eqref{schrodinger3} obeys
\beq \label{jit-last-condition2}
\limsup_{L \rightarrow \infty} \frac{\parallel u \parallel_L^2}{L^\beta}<\infty.
\eeq  
Then $\mu(A \cap \cdot)$ is continuous with respect to $(2-\beta)$-dimensional Hausdorff measure.
\end{proposition}

\begin{proposition}[\cite{jit-last}] \label{jit-last3}
Assume that $1 \leq a(j) \leq M$ for some $M>1$. If 
\beq \label{jit-last-condition3}
\liminf_{L \rightarrow \infty} \frac{\parallel u_{1,E} \parallel_L^2}{L^\alpha}=0
\eeq
for every $E$ in some Borel set $A$, then $\mu(A \cap \cdot)$ is singular with respect to $\alpha$-dimensional Hausdorff 
measure.
\end{proposition}

The next results we quote relate the properties of $\mu$ to the properties of the transfer matrices corresponding to $J$. 
These are the $2 \times 2$ matrices 
\beq \label{transfer-matrices1}
T_j(E)=S_j(E)S_{j-1}(E)\cdots S_1(E),
\eeq
where
\beq \label{transfer-matrices}
S_j(E)=\left( \begin{array}{cc}
\frac{E-b(j)}{a(j)} & -\frac{a(j-1)}{a(j)} \\
1 & 0 \end{array} \right). 
\eeq 

It isn't hard to see that 
\beq \label{transfer-solutions}
T_j(E)=\left( \begin{array}{cc}
u_{1,E}(j+1) & u_{2,E}(j+1) \\
u_{1,E}(j) & u_{2,E}(j) \end{array} \right) 
\eeq
so it is not surprising to find that the behavior of $T_j(E)$ is related to the behavior of the eigenfunctions.
For any $j_1,j_2$, we use the shorthand
\beq \label{t-j1-j2}
T_{j_1,j_2}(E) \equiv T_{j_1}(E)T_{j_2}(E)^{-1}.
\eeq

The following is a generalization of theorem 1.2 of \cite{last-simon} relating the behavior of $T_j(E)$ with the 
existence of absolutely continuous spectrum. 
\begin{proposition} \label{last-simon}
Let $m_j, l_j$ be arbitrary sequences of natural numbers and let
\beq \label{A1}
A_1= \left \{ E \mid \liminf_{j \rightarrow \infty} \frac{1}{a_{l_j}}\parallel T_{m_j,l_j}(E) \parallel < \infty \right \}.
\eeq
Then $A_1$ supports the a.c.\ part of $\mu$ in that $\mu_{\ac}(\bbR \setminus A_1)=0$.
\end{proposition}
\begin{proof}
Note that $\det(T_{l_j}(E))=\frac{1}{a_{l_j}}$ so that 
\beq \no
\parallel T_{l_j}(E)^{-1} \parallel = a_{l_j} \parallel T_{l_j}(E) \parallel.
\eeq
Thus, we have that 
\begin{align} \no
&\frac{1}{a_{l_j}}\parallel T_{m_j}(E) T_{l_j}(E)^{-1} \parallel \leq 
\frac{1}{a_{l_j}} \parallel T_{m_j}(E) \parallel \parallel T_{l_j}(E)^{-1} \parallel \notag \\
& =  \parallel T_{m_j}(E) \parallel \parallel T_{l_j}(E) \parallel. \notag
\end{align}
From here one proceeds exactly as in the proof of theorem 3.4D of \cite{last-simon} to get the conclusion:
We know that $u_{1,E}(j)$, viewed as functions of $E$, are orthonormal polynomials with respect to $d\mu$, and that 
$u_{2,E}(j)$ are orthonormal polynomials with respect to another measure - $d \tilde{\mu}$, such that the measure
$d\nu=\min(d\mu,d\tilde{\mu})$ is purely absolutely continuous and equivalent to the absolutely continuous part of $d\mu$ 
(see \cite{last-simon}). Thus, from the characterization \eqref{transfer-solutions}, we get that 
\beq \no
\int_{\bbR} \parallel T_{j}(E) \parallel^2 d\nu(E) \leq 4. 
\eeq 
The theorem now follows from an application of the Cauchy-Schwarz inequality and Fatou's lemma.
\end{proof}

In order to show absence of eigenvalues, the following idea of Simon-Stolz \cite{simon-stolz} is useful: 
\begin{proposition}[See \cite{simon-stolz}] \label{simon-stolz}
For a given $E \in \bbR$, if 
\beq \label{no-ev-condition}
\sum_{j, \det (T_j(E))=1} \parallel T_j(E) \parallel^{-2}=\infty,
\eeq
then there can be no eigenvalue at $E$.
\end{proposition}
\begin{proof}
The proposition follows from the observation that if 
\beq \no
\det (T_m(E))=1
\eeq 
and $u$ solves \eqref{schrodinger1} then 
\beq \no
|u(m+1)|^2+|u(m)|^2 \geq \frac{|u(1)|^2+|u(0)|^2}{\parallel T_m(E) \parallel^2}.
\eeq
\end{proof}

%%%%%%%%%%%%%%%%%%%%%%%%%%%%%%%%%%%%%%% REFERENCES %%%%%%%%%%%%%%%%%%%%%%%%%%%%%%%%%%%%%%%%%%%%%%%%%%%%%%%%%%%%%%%%%%%%%%

\end{document}